\documentclass[11pt,twoside,a4paper]{article}
\usepackage{amsfonts}
\usepackage{amssymb}
\usepackage{amsmath}

\begin{document}

\newcommand{\pa}{\partial}
\newcommand{\opa}{\overline\pa}
\newcommand{\ol}{\overline }

\numberwithin{equation}{section}

\newcommand\C{\mathbb{C}}  
\newcommand\R{\mathbb{R}}
\newcommand\Z{\mathbb{Z}}
\newcommand\N{\mathbb{N}}
\newcommand\PP{\mathbb{P}}

{\LARGE \centerline{Inflexible $CR$ submanifolds}}
\vspace{0.8cm}

\centerline{\textsc {Judith Brinkschulte\footnote{Universit\"at Leipzig, Mathematisches Institut, Augustusplatz 10, D-04109 Leipzig, Germany. 
E-mail: brinkschulte@math.uni-leipzig.de}}
 and C. Denson Hill
\footnote{Department of Mathematics, Stony Brook University, Stony Brook NY 11794, USA. E-mail: dhill@math.stonybrook.edu
 }}

\vspace{0.5cm}

\begin{abstract} In this paper we introduce the concept of {\it inflexible} $CR$ submanifolds. These are $CR$ submanifolds of some complex Euclidean space such that any compactly supported $CR$ deformation is again globally $CR$ embeddable into some complex Euclidean space. Our main result is that any  $2$-pseudoconcave quadratic $CR$ submanifold of type $(n,d)$ in $\C^{n+d}$ is inflexible.
\end{abstract}

\vspace{0.5cm}

\section{Introduction}

In this paper, we shall be interested in proving embedding results for compactly supported perturbations of embedded $CR$ manifolds.\\

Here an abstract $CR$ manifold of type $(n,d)$ is a triple $(M, HM, J)$, where $M$ is a smooth real manifold of dimension $2n+d$, $HM$ is a subbundle of rank $2n$ of the tangent bundle $TM$, and $J: HM \rightarrow HM$ is a smooth fiber preserving bundle isomorphism with $J^2= -\mathrm{Id}$. We also require that $J$ be formally integrable; i.e. that we have
$$\lbrack T^{0,1}M,T^{0,1}M\rbrack \subset T^{0,1}M$$
where 
$$ T^{0,1}M = \lbrace X+ iJX\mid X\in \Gamma(M,HM)\rbrace \subset \Gamma(M,\mathbb{C}TM),$$
with $\Gamma$ denoting smooth sections.

The $CR$ dimension of $M$ is $n\geq 1$ and the $CR$ codimension is $d\geq 1$.\\

A problem of great interest is to decide which $CR$ manifolds $M$ admit $CR$ embeddings into some complex Euclidean space. Namely, can one find a smooth embedding $\varphi$ of $M$ into $\mathbb{C}^N$ such that the induced $CR$ structure $\varphi_\ast(T^{0,1}M)$ on $\varphi(M)$ coincides with the $CR$ structure $T^{0,1}(\mathbb{C}^\N)\cap\mathbb{C}T(\varphi(M))$ from the ambient space $\mathbb{C}^N$.\\

Typically, examples of non-embeddable $CR$ structures arise as deformations of $CR$ submanifolds of some complex Euclidean space. For example, Rossi $\cite{R}$ constructed small real analytic deformations of the standard $CR$ structure on the 3-sphere $S^3$ in $\mathbb{C}^2$, and the resulting abstract $CR$ structures fail to $CR$ embed globally into $\mathbb{C}^2$. Also Nirenberg's famous local nonembeddability examples \cite{Ni} can be interpreted as small (local) deformations of the Heisenberg
structure on $\mathbb{H}^2\subset\mathbb{C}^2$. The examples by Nirenberg were later on extended to higher dimensions by Jacobwitz and Tr\`eves \cite{JT}.\\

However, there is something special about Nirenberg's three-dimensional examples: Since the formal integrability condition is always satisfied in this situation, one can easily modify the examples to obtain small (global) deformations of the Heisenberg structure $\mathbb{H}^2$. Moreover, these deformations are compactly supported (in the sense that the deformations coincide with the given Heisenberg structure outside a compact set). For the examples of Jacobowitz and Tr\`eves, it is not clear if this is  possible.\\

In fact, as soon as the $CR$ dimension is greater than one, the integrability conditions come into play, and they make it much more difficult to construct deformations. However, when $M$ is given as a $CR$ submanifold of some complex Eudlidean space, one can always obtain compact deformations of the $CR$ structure on $M$ by making a small compact geometric deformation of $M$ within the complex Euclidean space. We refer to this as "punching $M$". But it is not clear if there exists other compact deformations of the abstract $CR$ structure on $M$, which render $M$ no longer embeddable as a $CR$ submanifold of the complex Euclidean space, such as in Nirenberg's example. \\

Therefore in the present paper, we want to discuss the following problem: Suppose $f: (M,HM,J)\longrightarrow \C^{n+k}$ is a $CR$ embedding, and $(M^\prime, HM^\prime, J^\prime)$ is small, compactly supported $CR$ deformation of $(M,HM,J)$. Does it follow that it also admits a $CR$ embedding $f^\prime$ with $f^\prime$ close to $f$?.\\

An answer to this question clearly depends on the Levi-form of $M$, so let us now recall its intrinsic definition.\\

We denote by $H^o M=\lbrace \xi\in T^\ast M\mid < X,\xi>=0, \forall X\in H_{\pi(\xi)}M\rbrace$ the {\it characteristic conormal bundle} of $M$. Here $\pi: T M \longrightarrow M$ is the natural projection. To each $\xi\in H^o_p M\setminus \lbrace 0\rbrace$, we associate the Levi form at $\xi:$
$$\mathcal{L}_p(\xi, X) = \xi(\lbrack J\tilde X, \tilde X\rbrack )= d\tilde\xi(X,JX) \ \mathrm{for} \ X\in H_p M$$
which is Hermitian for the complex structure of $H_p M$ defined by $J$. Here $\tilde \xi$ is a section of $H^o M$ extending $\xi$ and $\tilde X$ a section of $HM$ extending $X$. \\

Following \cite{HN1} $M$ is called $q$-pseudoconcave, $0\leq q\leq\frac{n}{2}$ if for every $p\in M$ and every characteristic conormal direction $\xi\in H^o_p M\setminus \lbrace 0\rbrace$, the Levi form $\mathcal{L}_p(\xi, \cdot)$ has at least $q$ negative and $q$ positive eigenvalues.\\

{\bf Acknowledements.} The first author was supported by Deutsche Forschungsgemeinschaft (DFG, German Research Foundation, grant BR 3363/2-1).\\

\section{Definitions and statement of the main results}

Let $(M,HM,J)$ be $CR$ manifold of type $(n,d)$ globally $CR$ embedded into some complex Euclidean space. We say that $(M,HM,J)$ admits a {\it compactly supported $CR$ deformation} if there exists a family $(M_a, HM_a, J_a)_{a>0}$ of abstract $CR$ manifolds depending smoothly on a real parameter $a > 0$ and converging to $(M,HM,J)$ as $a$ tends to $0$ in the usual $\mathcal{C}^\infty$ topology; we also require that $(M_a, HM_a, J_a)= (M,HM,J)$ for every $a>0$ outside some compact $K$ of $M$ not depending on $a$.\\

We say that $(M,HM,J)$ is a {\it flexible} $CR$ submanifold if it admits a compactly supported $CR$ deformation $(M_a, HM_a, J_a)_{a>0}$ such that for every sufficiently small $a > 0$, the $CR$ structure $(M_a, HM_a, J_a)$ is not globally $CR$ embeddable into some complex Euclidean space. So, for example, the Heisenberg $CR$ structure $\mathbb{H}^2$ in $\C^2$ is flexible.\\

We say that $(M,HM,J)$ is an {\it inflexible} $CR$ submanifold if it is not flexible. That means that $(M,HM,J)$ is inflexible if and only if for every compactly supported $CR$ deformation $(M_a, HM_a, J_a)_{a>0}$ of $(M, HM, J)$, the $CR$ manifold $(M_a, HM_a, J_a)$ is globally $CR$ embeddable into some complex Euclidean space.\\

In other words, a flexible $CR$ submanifold admits a compactly supported $CR$ deformation that "pops out" of the space of globally $CR$ embeddable manifolds. On the other hand, for an inflexible $CR$ submanifold, any compactly supported $CR$ deformation stays in the space of globally $CR$ embeddable manifolds.\\

{\it Remark:} In the definitions above, we also allow compact deformations which are only defined for a sequence of $a$'s tending to zero.\\

Our main result is as follows:\\

\newtheorem{main}{Theorem}[section]
\begin{main}   \label{main}   \ \\
Let $M$ be a quadratic $CR$ submanifold of type $(n,d)$ in $\C^{n+d}$ that is $2$-pseudoconcave. Let $(M_a, HM_a, J_a)_{a>0}$ be a compactly supported $CR$ deformation of $(M,HM,J)$.
Then, given any smooth $CR$ function $f: (M,HM,J)\longrightarrow \C$, there is a $CR$ function $f_a: (M_a,HM_a,J_a)\longrightarrow \C$ as close to $f$ as we please, provided $a$ is sufficiently close to $0$.\\
Moreover, $f_a$ can be chosen to coincide with the given $f$ outside a compact of $M$. In particular, $(M_a, HM_a, J_a)$ is $CR$ embeddable into $\C^{n+d}$ for $a$ sufficiently close to $0$.
\end{main}

Here a {\it quadratic} $CR$ submanifold is a submanifold of $\C^{n+d}$ of the form
$$M=\lbrace z\in \C^{n+d}\mid \mathrm{Im} z_\ell = H_\ell(z_1,\ldots, z_n),\ n+1\leq \ell\leq n+d\rbrace,$$
where the $H_\ell$'s are quadratic hermitian forms on $\C^n$.\\

"$f_a$ as close to $f$ as we please" means that for any given $\ell\in\mathbb{N}$, any given compact $K$ of $M$ and arbitrary small $\varepsilon >0$, one can find a $CR$ function $f_a: (M_a,HM_a,J_a)\longrightarrow \C$ such that the $\mathcal{C}^\ell$ norm of $f-f_a$ on $K$  is less than $\varepsilon$.\\

In particular, Theorem \ref{main} implies that for a 2-pseudoconcave quadratic $CR$ submanifold, any compactly supported $CR$ deformation amounts to "punching $M$": any of the ambient complex coordinate functions is a $CR$ function on $M$. Our theorem yields that we can make arbitrarily small modifications of these coordinate functions inside a compact subset of $M$ to obtain global $CR$ coordinate functions on the deformed $CR$ manifolds.\\

The last statement of Theorem \ref{main} combined with the definition of "inflexible" immediately gives the following

\newtheorem{corr}[main]{Corollary}
\begin{corr}   \label{corr}   \ \\
Let $M$ be 2-pseudoconcave quadratic $CR$ submanifold of type $(n,d)$ in $\C^{n+d}$. Then $M$ is inflexible.
\end{corr}

\section{A first example}

The idea  of the proof of Theorem \ref{main} is as follows: For a given $CR$ function $f$ on $M$ we want to find a $CR$ function on $M_a$ which is very close to the given $f$. Therefore we want to solve the Cauchy-Riemann equations  $\opa_{M_a}u = \opa_{M_a}f $ with $u$ having compact support and the $\mathcal{C}^k$-norms of $u$ being controlled by some $\mathcal{C}^l$-norms of $\opa_{M_a}$ (uniformly with respect to $a$). \\

In this section, we  will explicitly carry out the proof of our main result \ref{main} in all details for the easiest example of a 2-pseudoconcave $CR$ manifold. Namely let $M\subset \C^5$ be the real hypersurface defined by
\begin{equation}  \label{defM}
M = \lbrace (z_1,z_2,z_3, z_4, x + iy) \mid y = \vert z_1\vert^2 + \vert z_2\vert^2 - \vert z_3\vert^2 - \vert z_4\vert^2 \rbrace.
\end{equation}
Then $M$ is a 2-pseudoconcave $CR$ manifold of type $(4,1)$. To abbreviate notations, we also define $z= (z_1, z_2, z_3, z_4)$ and $\vert z\vert^2 = \vert z_1\vert^2 + \vert z_2\vert^2 + \vert z_3\vert^2 + \vert z_4\vert^2$.  A straightforward computation shows that $T^{0,1}M$ is spanned by
$$\ol L_j = \frac{\pa}{\pa\ol z_j} - i \epsilon_j z_j\frac{\pa}{\pa x}, \ j =1,2,3,4,$$
where $\epsilon_1 =\epsilon_2 =-1$ and $\epsilon_3 = \epsilon_4 = 1$. \\
For $u = \sum_{1}^{4} u_j d\ol z_j\in \mathcal{C}^\infty_{0,1}(M)$ we have
$$\opa_M u = \sum_{j,k=1}^4 \ol L_k(u_j) d\ol z_k\wedge d\ol z_j.$$
Next, we consider the volume element 
$$dV = (\frac{i}{2})^4 e^{\vert z\vert^2} \bigwedge_{j=1}^4 dz_j\wedge d\ol z_j \wedge dx = \frac{1}{16} e^{\vert z\vert^2} \bigwedge_{j=1}^4 dz_j\wedge d\ol z_j \wedge dx$$
 on $M$, and we denote by $\Vert\ \Vert$ the $L^2$-norm of $(0,q)$-forms on $M$ with respect to this volume element, where the pointwise norms of $(0,q)$-forms on $M$ is the one induced by the standard euclidean metric on $\mathbb{C}^5$. The corresponding $L^2$-spaces will be denoted by $L^2_{0,q}(M, \vert z\vert^2)$. Then $\opa^\ast_{M}$, the formal adjoint of $\opa_M$ with respect to $\Vert\ \Vert$ can be computed as follows:
$$\opa^\ast_{M}u = - e^{-\vert z\vert^2} \sum_{j=1}^4 L_j(u_j e^{\vert z\vert^2})$$
for $u\in \mathcal{D}^{0,1}(M)$. \\

First we will prove the following $L^2$ estimates on $M$:\\

\newtheorem{L2}{Lemma}[section]  
\begin{L2}  \label{L2} 
\ \\ Let $M$ be defined as in (\ref{defM}).\begin{enumerate}
\item For all $u\in L^2_{0,1}(M, \vert z\vert^2) \cap \mathrm{Dom}(\opa_M)\cap\mathrm{Dom}(\opa_M^\ast)$ we have
 \begin{equation}  \label{vanishing}
2 \Vert u\Vert^2 \leq \Vert \opa_M u\Vert^2 + \Vert\opa^\ast_{M} u\Vert^2.
\end{equation}
\item For all $u\in L^2_{0,0}(M, \vert z\vert^2) \cap \mathrm{Dom}(\opa_M)$ we have
\begin{equation}  \label{noCR}
4 \Vert u\Vert^2 \leq \Vert \opa_M u\Vert^2
\end{equation}
\end{enumerate}
\end{L2}

{\it Proof.} Throughout the proof of this Lemma, we identify $M$ with $\mathbb{C}^4\times\mathbb{R}$.
We will begin by showing how to reduce the proof of (\ref{vanishing}) to an estimate for an easier differential operator. Therefore we introduce the partial Fourier transform with respect to the variable $x$:
$$\tilde u(z,\xi)= \int e^{-i \langle x,\xi \rangle} u(z,x) dx$$
(for differential forms, this partial Fourier transform is defined componentwise).\\

Now an easy computation shows that for $u\in\mathcal{D}^{0,1}(M)$ we have
\begin{eqnarray*}
\widetilde{\opa_M u}(z,\xi) & = & \sum_{j,k=1}^4 \widetilde{\ol L_k(u_j)}(z,\xi)d\ol z_k\wedge d\ol z_j\\
& = & \sum_{j,k=1}^4 \big(\frac{\pa}{\pa\ol z_k}u_j - i\epsilon_k z_k\frac{\pa}{\pa x}u_j\big)\widetilde{\ \ \ }(z,\xi)d\ol z_k\wedge d\ol z_j\\
& = & \sum_{j,k=1}^4 \big( \frac{\pa}{\pa\ol z_k}\tilde{u}_j(z,\xi) + \epsilon_k z_k \xi \tilde{u}_j(z,\xi) \big) d\ol z_k\wedge d\ol z_j\\
& = & \opa_{(z)} \tilde u(z,\xi),
\end{eqnarray*}
where $\opa_{(z)}$ is defined by
$$\opa_{(z)}v(z,\xi)= \sum_{j,k=1}^4 \opa_k v_j  d\ol z_k\wedge d\ol z_j.$$
Here $\opa_k v_j= \frac{\pa}{\pa\ol z_k}v_j +\epsilon_k z_k \xi v_j$ is of order 0 in $\xi$. Similarly, we get 
\begin{eqnarray*}
\widetilde{\opa^\ast_{M}u}(z,\xi) & = & - (\sum_{j=1}^4 \widetilde{ L_j u_j +  \ol z_j u_j})(z,\xi)\\
& = & \delta_{(z)} \tilde{u} (z,\xi),
\end{eqnarray*}
where
$$\delta_{(z)}v(z,\xi)= \sum_{j=1}^4(\delta_j v_j) (z,\xi) $$
with $\delta_j v_j = -\frac{\pa}{\pa z_j} v_j + \epsilon_j \ol z_j \xi v_j - \ol z_j v_j$.
Note that also $\delta_j$ is of order 0 in $\xi$.\\

Now, as in \cite{H1} we compute
\begin{eqnarray} \label{first}
\vert \opa_{(z)} v\vert^2 & = & \vert \sum_{j,k=1}^4\opa_k v_j d\ol z_k\wedge d\ol z_j\vert^2\nonumber\\
& = & \frac{1}{2} \sum_{j,k=1}^4 \vert \opa_k v_j - \opa_j v_k\vert^2 \nonumber\\
& = &  \sum_{j,k=1}^4 \vert \opa_j v_k\vert^2 -  \sum_{j,k=1}^4 \opa_k v_j\ol{\opa_j v_k}
\end{eqnarray}

Also we have
\begin{eqnarray}  \label{second}
\vert \delta_{(z)}v\vert^2 & = & \vert \sum_{j=1}^4\delta_j v_j\vert^2 = \sum_{j,k=1}^4 \delta_j v_j \ol{\delta_k v_k} \nonumber\\
& = & \sum_{j=1}^4 \vert \delta_j v_j\vert^2 + \sum_{j\not= k}\delta_j v_j \ol{\delta_k v_k}
\end{eqnarray}

Summing up (\ref{first}) and (\ref{second}) we obtain

$$\int_{\mathbb{C}^4} \big( \vert\opa_{(z)}v\vert^2 + \vert \delta_{(z)}v\vert^2\big)\exp{(\vert z\vert^2)} (\frac{i}{2})^4\bigwedge_{j=1}^{4}dz_j\wedge d\ol z_j =$$
 $$\sum_{j=1}^4 \Vert\delta_j v_j\Vert^2_z + \sum_{j\not=k} \Vert\opa_j v_k\Vert^2_z + \sum_{j\not=k} \ll \lbrack \opa_k,\delta_j\rbrack v_j, v_k\gg_z .$$
Here we have used that $\opa_k$ and $\delta_k$ are adjoint operators. To abbreviate notations, we have introduced $\Vert\ \Vert_z$ to denote partial integration with respect to the $z=(z_1,z_2,z_3,z_4)$ variables:
$$\Vert v\Vert_z^2 = \int_{z\in \mathbb{C}^4} \vert v(z,\xi)\vert^2 \exp{(\vert z\vert^2)} (\frac{i}{2})^4\bigwedge_{j=1}^{4}dz_j\wedge d\ol z_j.$$

 Since $\lbrack \opa_k,\delta_j\rbrack = 0$ for $j\not= k$ we obtain
\begin{equation}  \label{apriori1}
\Vert\opa_{(z)}v\Vert^2_z + \Vert \delta_{(z)}v\Vert^2_z = 
\sum_{j\not= k} \Vert\opa_j v_k\Vert^2_z + \sum_{j=1}^4\Vert \delta_j v_j\Vert^2_z 
\end{equation}

Also, a straightforward computation shows that
\begin{equation}  \label{commutator}
\lbrack\opa_j,\delta_j\rbrack = -1 +2\epsilon_j \xi
\end{equation}

This will be used to show that for each fixed $k\in \lbrace 1,2,3,4\rbrace$ we have
\begin{equation}  \label{apriori2}
\sum_{\underset{j\not= k}{j=1} }^4 \Vert \opa_j v_k\Vert^2_z \geq 2 \Vert v_k\Vert^2_z
\end{equation}

Assume e.g. k=4.  From (\ref{commutator})  we then obtain
\begin{equation}   \label{identity}
\Vert \delta_j v_4\Vert^2_z - \Vert \opa_j v_4\Vert^2_z = (-1+2\epsilon_j \xi) \Vert v_4\Vert^2_z.
\end{equation}
It follows that
\begin{eqnarray*}
\sum_{j=1}^3 \Vert \opa_j v_4\Vert^2_z & \geq & \sum_{j=1,3} \Vert \opa_j v_4\Vert^2_z\\
& = & \sum_{j=1,3} \Vert \delta_j v_4\Vert^2_z + 
\sum_{j=1,3} (1- 2\epsilon_j \xi) \Vert v_4\Vert^2_z\\
& \geq & 2\Vert  v_4\Vert^2_z - 2\xi(-1+1) \Vert v_4\Vert^2_z\\
& = & 2 \Vert v_4\Vert^2_z,
\end{eqnarray*}
which proves (\ref{apriori2}) for $k=4$. The remaining cases are similar.\\

Combining (\ref{apriori1}) and (\ref{apriori2}) we have proved that
$$2\Vert v\Vert_z^2(\xi) \leq \Vert\opa_{(z)} v\Vert_z^2(\xi) + \Vert \delta_{(z)}v\Vert_z^2(\xi)$$
for every fixed $\xi\in\mathbb{R}$. Setting $v= \tilde{u}$ and integrating this inequality with respect to $\xi$ we obtain from the definition of the operators $\opa_{(z)}$ and $\delta_{(z)}$
$$2\Vert \tilde{u}\Vert^2 \leq \Vert \widetilde{\opa_M u}\Vert^2 + \Vert \widetilde{\opa_M^\ast u}\Vert^2$$
for all $u\in\mathcal{D}^{0,1}(M)$.\\

The Plancherel theorem permits to conclude that
$$2\int_M \vert u\vert^2 dV \leq \int_M( \vert\opa_M u\vert^2 + \vert\opa^\ast_{M}\vert^2) dV$$
for $u\in\mathcal{D}^{0,1}(M)$. Obviously, the restriction of the standard euclidean metric to $M$ is complete, therefore the above estimate extends to all $u\in L^2_{0,1}(M, \vert z\vert^2) \cap \mathrm{Dom}(\opa_M)\cap\mathrm{Dom}(\opa_M^\ast)$, which proves the first statement of the Lemma.\\

The proof of (\ref{noCR}) is similar. Indeed, using the partial Fourier transform, the proof of (\ref{noCR}) is again reduced to the estimate of $\sum_{j=1}^4 \Vert \opa_j v\Vert^2_z$, where $\opa_j$ is defined as before. But using (\ref{identity}) we get
\begin{eqnarray*}
\sum_{j=1}^4 \Vert \opa_j v\Vert^2_z & = & \sum_{j=1}^4 \Vert \delta_j v\Vert^2_z + 
\sum_{j=1}^4 (1- 2\epsilon_j \xi)\Vert v\Vert^2_z\\
& \geq & 4 \Vert  v\Vert^2_z - 2\xi(-1-1+1+1) \Vert v\Vert^2_z\\
& = & 4 \Vert v\Vert^2_z,
\end{eqnarray*}
This completes the proof of the Lemma by the same arguments as before.\hfill $\square$\\

Next, we use again that $M$ is 2-pseudoconcave (this condition is clearly stable under small perturbations). This implies that we have a  uniform subelliptic estimate in degree $(0,1)$ (see \cite{FK}):

 For every compact $K$ of $M$, there exists a constant $C_K > 0$ independent of $a$ such that
\begin{equation}  \label{subelliptic}
\Vert u\Vert^2_{\frac{1}{2}} \leq C_K (\Vert \opa_{M_a} u\Vert^2 + \Vert\opa^\ast_{M_a} u\Vert^2 + \Vert u\Vert^2)
\end{equation}
for all $u\in \mathcal{D}_K^{0,1}(M_a)$.\\

Combining Lemma \ref{L2} and (\ref{subelliptic}), we can establish an $L^2$ a priori estimate in degree $(0,1)$, which is uniform with respect to $a$ (in the sense that the  constant involved does not depend on $a$).

\newtheorem{firstlemma}[L2]{Lemma}
\begin{firstlemma}   \label{firstlemma}   \ \\
There is  $a_0> 0$ and a constant $C > 0$ such that
$$\Vert u \Vert^2\leq C (\Vert \opa_{M_a} u\Vert^2 + \Vert \opa^\ast_{M_a}u\Vert^2)$$
for all $u\in L^2_{0,1}(M_a,\vert z\vert^2)$, $a < a_0$.
\end{firstlemma}

{\it Proof.} Following \cite{N},
assume by contradiction that there is a sequence $\lbrace u_{a_\nu}\rbrace\in L^2_{0,1}(M_{a_{\nu}},\vert z\vert^2)\cap\mathrm{Dom}(\opa_{M_{a_{\nu}}})\cap \mathrm{Dom}(\opa^\ast_{M_{a_{\nu}}}) $, $a_\nu \rightarrow 0$, such that 
\begin{equation}  \label{1}
\Vert u_{a_\nu}\Vert = 1,
\end{equation}
whereas
\begin{equation}  \label{2}
\Vert \opa_{M_{a_\nu}} u_{a_\nu}\Vert + \Vert \opa^\ast_{M_{a_\nu}} u_{a_\nu}\Vert < a_\nu.
\end{equation}

We now want to show that $\lbrace u_{a_\nu}\rbrace$ is a Cauchy sequence. \\

Remember that $M_{a_\nu}= M$ outside $K$. We now choose a slightly larger compact $K_1$ containing $K$ in its interior, and a smooth cut-off function $\chi$ such that $\chi\equiv 1$ outside $K_1$ and $\chi\equiv 0$ in a neighborhood of $K$. Since $\opa_{M_{a_\nu}}$, $\opa^\ast_{M_{a_\nu}}$ coincide with $\opa_M$, $\opa^\ast_{M}$ outside $K$, we obtain from (\ref{vanishing})
$$2\Vert \chi u\Vert^2 \leq  \Vert \opa_M (\chi u)\Vert^2 + \Vert\opa^\ast_{M}(\chi u)\Vert^2 $$
for all $u\in L^2_{0,1}(M_a,\vert z\vert^2)$, which implies
\begin{equation}  \label{outsideK}
\Vert \chi u\Vert^2 \leq C^\prime (\Vert \opa_M  u\Vert^2 + \Vert\opa^\ast_{M}u\Vert^2 + \int_{K_1\setminus K} \vert u\vert^2 dV )
\end{equation}
for some constant $C^\prime > 0$.\\

On the other hand, let $\eta$ be a smooth cut-off function so that $\eta\equiv 1$ in a neighborhood of $K_1$. Then $\Vert \eta u_{a_\nu} \Vert_{\frac{1}{2}}$ is bounded by (\ref{subelliptic}), so the generalized Rellich lemma implies that the sequence $\lbrace u_{a_\nu}\rbrace$ restricted to $K_1$ is precompact in $L^2_{0,1}(K_1)$. Thus it is no loss of generality to asume that the restriction of $\lbrace u_{a_\nu}\rbrace$ to $K_1$ is a Cauchy sequence. But this combined with (\ref{outsideK}) implies that $\lbrace u_{a_\nu}\rbrace$ is a Cauchy sequence in $L^2_{0,1}(M,\vert z\vert^2)$.\\

 Denote by $u_0$ the limit of this sequence. From (\ref{2}) it follows that $\opa_M u_0$ and $\opa^\ast_{M}u_0$, defined in the distribution sense, both vanish. But from (\ref{1}) it also follows that $\Vert u_0\Vert = 1$. This contradicts (\ref{vanishing}) and therefore completes the proof of the lemma. \hfill$\square$\\

{\it Proof of theorem \ref{main} for $M$ as above.}

Let $f$ be given. Then $\opa_{M_a}f$ has compact support and  tends to zero when $a$ tends to zero. It is well known (see e.g. \cite{H2}) that the a priori estimate (\ref{vanishing}) implies that we can solve the equation $\opa_{M_a} u_a = \opa_{M_a} f$ with $\Vert u_a\Vert \leq C\Vert\opa_{M_a}u_a\Vert$. Hence $u_a$ is as small as we wish in $L^2(M_a, \vert z\vert^2)$, provided $a$ is small enough. It is well-known that the subelliptic estimate (\ref{subelliptic}) implies also the following:  Suppose given a compact $K\subset M_a$  and two smooth real functions $\zeta,\ \zeta_1$ with 
$\mathrm{supp}\zeta \subset\mathrm{supp}\zeta_1\subset K$ and $\zeta_1 =1$ on $\mathrm{supp}\zeta$, then
for any integer $m\in\mathbb{N}$ there exists a constant $C_{K,m}$ such that
$$\Vert \zeta u\Vert^2_{m+\varepsilon} \leq C_{K,m} (\Vert \zeta_1\opa_{M_a}u\Vert^2_m + \Vert \zeta_1\opa^\ast_{M_a} u\Vert^2_m + \Vert \zeta_1 u\Vert^2)$$
Here $\Vert\ \Vert_m$ denotes the Sobolev norm of order $m$. But then, choosing the minimal solution satisfying $\opa^\ast_{M_a}u=0$, also
  the $\mathcal{C}^\ell$-norm of $u_a$ over a given compact $K\subset M_a$ can be controlled by some $\mathcal{C}^m$-norm of $\opa_{M_a}u_a = f$, and hence made small when letting $a$ tend to zero.
Setting $f_a = f- u_a$  proves the first statement.\\

Moreover,  $u_a$ has compact support: Since the $CR$ structures of $M$ and $M_a$ coincide outside a compact set, and $u_a$ solves the equation $\opa_{M_a}u_a= \opa_{M_a}f$, $u_a$ is a $CR$ function on $M$ outside some compact set $K$. It is no loss of generality to assume that $M\setminus K$ is connected. But then, since the Hartogs phenomenon for $CR$ functions holds in 2-pseudoconcave $CR$ manifolds \cite{LT}, the restriction of $u_a$ to $M\setminus K$ extends to a $CR$ function $\tilde u_a$ on $M$. Since $u_a$ belongs to $L^2_{0,0}(M,\vert z\vert^2)$, the same is true for $\tilde u_a$. But then (\ref{noCR}) implies $\tilde u_a\equiv 0$. Hence $u_a$ vanishes on $M\setminus K$. \hfill$\square$\\

\section{The general case}

In this section we will explain the proof of Theorem \ref{main} for a general 2-pseudoconcave quadratic $CR$ submanifold $M$ of type $(n,d)$ given by
$$M=\lbrace z\in \C^{n+d}\mid \mathrm{Im} z_\ell = \sum_{i,j =1}^n h^\ell_{ij} z_i\ol z_j,\ n+1\leq \ell\leq n+d\rbrace.$$

In this case, $T^{1,0}M$ is spanned by
$$ L_j = \frac{\pa}{\pa z_j} + i\sum_{\ell = n+1}^{n+d}\sum_{k=1}^n h^\ell_{jk}\ol z_k \frac{\pa}{\pa x_\ell}      \quad j=1,\ldots,n,$$
and $T^{0,1}M$ is spanned by
$$ \ol L_j = \frac{\pa}{\pa\ol z_j} - i\sum_{\ell = n+1}^{n+d}\sum_{k=1}^n h^\ell_{kj} z_k \frac{\pa}{\pa x_\ell}      \quad j=1,\ldots,n.$$

First we show that the analogue of Lemma \ref{L2} still holds true, i.e. we have the following

\newtheorem{L2general}{Lemma}[section]  
\begin{L2general}  \label{L2general} 
\ \\ Let $M$ be a $2$-pseudoconcave quadratic $CR$ submanifold.\begin{enumerate}
\item For all $u\in L^2_{0,1}(M, \vert z\vert^2) \cap \mathrm{Dom}(\opa_M)\cap\mathrm{Dom}(\opa_M^\ast)$ we have
 \begin{equation}  \label{vanishinggen}
\Vert u\Vert^2 \leq \Vert \opa_M u\Vert^2 + \Vert\opa^\ast_{M} u\Vert^2.
\end{equation}
\item For all $u\in L^2_{0,0}(M, \vert z\vert^2) \cap \mathrm{Dom}(\opa_M)$ we have
\begin{equation}  \label{noCRgen}
 \Vert u\Vert^2 \leq \Vert \opa_M u\Vert^2
\end{equation}
\end{enumerate}
\end{L2general}

{\it Proof of Lemma \ref{L2general}.} We show how the proof of Lemma \ref{L2} generalizes to this more general setting. In fact, we again use the partial Fourier transform with respect to the variables $(x_{n+1},\ldots,x_{n+d})$.
For a fixed $\xi\in\mathbb{R}^d$, we define the hermitian matrix 
$$h^\xi = \sum_{\ell= n+1}^d H_\ell\xi_\ell,\quad\quad\mathrm{i.e.}\ h^\xi_{jk}= \sum_{\ell= n+1}^d h^\ell_{jk}\xi_\ell.$$
After possibly making a unitary change of coordinates in the variables $(z_1,\ldots,z_n)$, we may assume that $h^\xi$ is diagonal with diagonal entries $h^\xi_{jj} = \lambda_j$ with $\lambda_1\leq \ldots\leq \lambda_n$.\\

Then, as in the proof of Lemma \ref{L2} we compute $\widetilde{\opa_M u}(z,\xi) = \opa_{(z)}\tilde u (z,\xi)$ with 
$$\opa_{(z)}v(z,\xi) = \sum_{k,s=1}^n\opa_k v_s d\ol z_k\wedge d\ol z_s,$$
where
\begin{eqnarray*}
\opa_k v_s & = & \frac{\pa}{\pa\ol z_k}v_s + \sum_{\ell=n+1}^{n+d}\sum_{m=1}^n h^\ell_{mk}z_m\xi_\ell v_s \\
& = & \frac{\pa}{\pa\ol z_k}v_s + \sum_{m=1}^n h^\xi_{mk}z_m v_s\\
& = & \frac{\pa}{\pa\ol z_k}v_s + \lambda_k z_k v_s
\end{eqnarray*}

Similarly we get
\begin{eqnarray*}
\widetilde{\opa^\ast_{M}u}(z,\xi) 
& = & \delta_{(z)} \tilde{u} (z,\xi),
\end{eqnarray*}
where
$$\delta_{(z)}v(z,\xi)= \sum_{j=1}^n(\delta_j v_j) (z,\xi) $$
with 
\begin{eqnarray*}
\delta_j v_j & = & -\frac{\pa}{\pa z_j} v_j + \sum_{\ell = n+1}^{n+d}\sum_{k=1}^n h^\ell_{jk}\ol z_k\xi_\ell v_j - \ol z_j v_j\\
& = & -\frac{\pa}{\pa z_j} v_j + \sum_{k=1}^n h^\ell_{jk}\ol z_k v_j - \ol z_j v_j\\
& = & -\frac{\pa}{\pa z_j} v_j + \sum_{k=1}^n \lambda_j \ol z_j v_j - \ol z_j v_j .
\end{eqnarray*}

The commutator of $\opa_k$ and $\delta_j$ can be computed as
\begin{equation}  \label{comm}
\lbrack\opa_k,\delta_j\rbrack = (-1+2\lambda_j) \delta_{j,k} ,
\end{equation}
where $\delta_{j,k}$ denotes the Kronecker symbol.\\

Therefore, as in the proof of Lemma \ref{L2}, one obtains for $v\in\mathcal{D}^{0,1}(\mathbb{C}^n\times\mathbb{R}^d)$:
\begin{eqnarray}  \label{firstest}
\Vert\opa_{(z)}v\Vert^2_z + \Vert \delta_{(z)}v\Vert^2_z  & = &
\sum_{j\not= k} \Vert\opa_j v_k\Vert^2_z + \sum_{j=1}^n\Vert \delta_j v_j\Vert^2_z \nonumber\\
& \geq & \sum_{j\not= k} \Vert\opa_j v_k\Vert^2_z.
\end{eqnarray}

Now we fix $k\in\lbrace 1,\ldots,n\rbrace$. Since $M$ is $2$-pseudoconcave, the hermitian matrix $h^\xi$ has at least $2$ negative and $2$ positive eigenvalues. But this implies that there exist indices $r,s \not=k$ such that $\lambda_r < 0$ and $\lambda_s > 0$. We now define real numbers $a_j\in\lbrack 0,1\rbrack $ by
\begin{eqnarray*}
a_j & = & 0,\ j\not= r,s \\
a_r & = & \frac{\lambda_s}{\lambda_s -\lambda_r} \\
a_s & = & \frac{-\lambda_r}{\lambda_s -\lambda_r}
\end{eqnarray*}
Note that by definition of $a_j$ we have $\sum_{j=1}^n a_j = 1$ and
$\sum_{j=1}^n a_j \lambda_j =0$. 
But then, using (\ref{comm}) we obtain
\begin{eqnarray*}
\sum_{\underset{j\not=k}{j=1}}^n \Vert\opa_j v_k\Vert^2_z & \geq & \sum_{\underset{j\not=k}{j=1}}^n a_j \Vert\opa_j v_k\Vert^2_z\\
& = & \sum_{j=1}^n a_j \Vert\delta_j v_k\Vert^2_z + \sum_{j=1}^n (1-2\lambda_j) a_j \Vert v_k\Vert^2_z \\
& \geq & \sum_{j=1}^n a_j \Vert v_k\Vert^2_z -2\sum_{j=1}^n \lambda_j a_j \Vert v_k\Vert^2_z\\
& \geq & \Vert v_k\Vert^2_z.
\end{eqnarray*}
From (\ref{firstest}) we therefore obtain
$$\Vert\opa_{(z)}v\Vert^2_z + \Vert \delta_{(z)}v\Vert^2_z \geq \Vert v\Vert^2_z.$$
By reasoning as in the proof of Lemma \ref{L2} we may therefore conclude that (\ref{vanishinggen}) holds.\\

Likewise, for the proof of (\ref{noCRgen}), we
define real numbers $c_j\in\lbrack 0,1\rbrack $ by
\begin{eqnarray*}
c_j & = & 0,\ j\not= 1,n \\
c_1 & = & \frac{\lambda_n}{\lambda_n -\lambda_1} \\
c_n & = & \frac{-\lambda_1}{\lambda_n -\lambda_1}
\end{eqnarray*}
Then we have $\sum_{j=1}^nc_j =1$ and $\sum_{j=1}^nc_j\lambda_j =0$. Therefore (\ref{comm}) implies
\begin{eqnarray*}
\sum_{j=1}^n \Vert \opa_j v\Vert^2_z & \geq  & 
\sum_{j=1}^n c_j \Vert \opa_j v\Vert^2_z\\
& = & 
\sum_{j=1}^n c_j\Vert \delta_j v\Vert^2_z + 
\sum_{j=1}^n c_j (1- 2 \lambda_j )\Vert v\Vert^2_z\\
& \geq & \sum_{j=1}^n c_j\Vert  v\Vert^2_z -2
\sum_{j=1}^n c_j  \lambda_j \Vert v\Vert^2_z\\
& = &  \Vert v\Vert^2_z,
\end{eqnarray*}
This completes the proof of (\ref{noCRgen}) by the same arguments as in the proof of Lemma \ref{L2}.\hfill $\square$\\

{\it Remark:} The proof of this Lemma is essentially contained in \cite{N} with constants depending on the Levi form of $M$. Here we have shown that one can take the same constant $1$ for every $2$-pseudoconcave quadratic $CR$ submanifold $M$.\\

The second essential ingredient for the proof of Theorem \ref{main} in the general case is the subelliptic estimate proved for 2-pseudoconcave $CR$ manifolds of arbitrary codimension $d$ in \cite{HN1}:
There exists  $\varepsilon > 0$ such that
for every compact $K$ of $M$, there exists a constant $C_K > 0$ independent of $a$ such that
\begin{equation}  \label{subellipticgen}
\Vert u\Vert^2_{\varepsilon} \leq C_K (\Vert \opa_{M_a} u\Vert^2 + \Vert\opa^\ast_{M_a} u\Vert^2 + \Vert u\Vert^2) 
\end{equation}
for all $u\in \mathcal{D}_K^{0,1}(M_a)$.
This subelliptic estimate replaces (\ref{subelliptic}) in the general situation.\\

Using (\ref{vanishinggen}) and (\ref{subellipticgen}), one can prove the uniform $L^2$ a priori estimate for $\opa_{M_a}$ as  stated in Lemma \ref{firstlemma}. The proof is the same. But this, together with (\ref{noCRgen}) completes the proof of Theorem \ref{main} as in section 3.

\end{document}